\documentclass[12pt]{article}

\usepackage{latexsym}

\newtheorem{theorem}{Theorem}[section]

\newtheorem{example}[theorem]{Example}
\newtheorem{definition}[theorem]{Definition}

\newtheorem{proposition}[theorem]{Proposition}

\newtheorem{lemma}[theorem]{Lemma}

\newenvironment{proof}{\medskip\noindent{\it Proof.\ }}{\hfill \mbox{$\Box$}\medskip}

\begin{document}

\def\eqnsep{60pt}
\def\enumsep{30pt}

\title{Restricted Permutations Related to Fibonacci Numbers and $k$-Generalized Fibonacci Numbers\footnote{MR Subject Classification:  05A15}}

\author{Eric S. Egge\\
\\
Department of Mathematics \\
Gettysburg College \\
300 North Washington Street \\
Gettysburg, PA  17325  USA \\
\\
eggee@member.ams.org}

\date{}

\maketitle

\begin{abstract}
A permutation $\pi \in S_n$ is said to {\it avoid} a permutation $\sigma \in S_k$ whenever $\pi$ contains no subsequence with all of the same pairwise comparisons as $\sigma$.
In 1985 Simion and Schmidt showed that the number of permutations in $S_n$ which avoid 123, 132, and 213 is the Fibonacci number $F_{n+1}$.
In this paper we generalize this result in two ways.
We first show that the number of permutations which avoid 132, 213, and $12\ldots k$ is the $k-1$-generalized Fibonacci number $F_{k-1,n+1}$.
We then show that the number of permutations which avoid 123, 132, and $k-1\ k-2\ \ldots 2 1 k$ is also the $k-1$-generalized Fibonacci number $F_{k-1,n+1}$.
We go on to show that the number of permutations in $S_n$ which avoid 132, $k\ k-1\ \ldots 4 2 1 3$, and 2341 is given by a polynomial plus a linear combination of two Fibonacci numbers.
We give explicit enumerations for $k \le 6$.
We begin to generalize this result by showing that the number of permutations in $S_n$ which avoid 132, 213, and $23\ldots k1$ is $\sum_{i=1}^n F_{k-2,i}$.
We conclude with several conjectures and open problems.

\medskip

\end{abstract}

\section{Introduction}

Let $S_n$ denote the set of permutations of $[n] = \{1, \ldots, n\}$, written in one-line notation, and suppose $\pi \in S_n$ and $\sigma \in S_k$.
We say $\pi$ {\it avoids} $\sigma$ whenever $\pi$ contains no subsequence with all of the same pairwise comparisons as $\sigma$.
For example, the permutation 23481756 avoids 4213, but it has 4156 as a subsequence so it does not avoid 2134.
If $\pi$ avoids $\sigma$ then $\sigma$ is sometimes called a {\it pattern} or a {\it forbidden subsequence} and $\pi$ is sometimes called a {\it restricted permutation} or a {\it pattern-avoiding permutation}.

One important and often difficult problem in the study of restricted permutations is the enumeration problem:  given a set $R$ of permutations, enumerate the set $S_n(R)$ consisting of those permutations in $S_n$ which avoid every element of $R$.
The earliest solution to an instance of this problem seems to be MacMahon's enumeration of $S_n(123)$, which is implicit in chapter V of \cite{MacMahon}.
The first explicit solution seems to be Hammersley's enumeration of $S_n(321)$ in \cite{Hammersley}.
In \cite{KnuthVol1} and \cite{KnuthVol3} Knuth shows that for any pattern $\sigma$ of length three, $S_n(\sigma) = C_n$, the $n$th Catalan number.
Other authors considered restricted permutations in the 1970s and early 1980s (see, for instance, \cite{Rogers}, \cite{Rotem1}, and \cite{Rotem2}) but the first systematic study was not undertaken until 1985, when Simion and Schmidt \cite{SimionSchmidt} solved the enumeration problem for every subset of $S_3$.
More recent work on various instances of the enumeration problem may be found in \cite{Atkinson}, \cite{BonaExact}, \cite{BonaSmooth}, \cite{Kremer1}, \cite{KremerShiu}, \cite{Mansour33k}, \cite{MansourVainshtein}, \cite{MansourVainshtein2}, \cite{Stankova1}, \cite{Stankova2}, \cite{StankovaWestHex}, \cite{WestCatalanSchroder}, and \cite{WestGenTrees}.

In this paper we are concerned with instances of the enumeration problem whose solutions involve the Fibonacci numbers or the $k$-generalized Fibonacci numbers.
The earliest example of such a result is Simion and Schmidt's proof \cite[Prop. 15]{SimionSchmidt} that 
\begin{equation}
\label{eqn:SimionSchmidtintro}
S_n(123, 132, 213) = F_{n+1} \hspace{\eqnsep} (n \ge 1).
\end{equation}
Somewhat later West showed \cite{WestGenTrees} that for many sets $R$ consisting of one pattern of length three and one of length four, $|S_n(R)| = F_{2n-1}$.
More recently Mansour showed \cite[Ex. 1, 3]{Mansour33k} that
\begin{equation}
\label{eqn:MansourintroFib1}
|S_n(123, 132, 3241)| = F_{n+2} - 1 \hspace{\eqnsep} (n \ge 1),
\end{equation}
\begin{equation}
\label{eqn:MansourintroFib2}
|S_n(132, 213, 2341)| = F_{n+2} - 1 \hspace{\eqnsep} (n \ge 1),
\end{equation}
and
\begin{equation}
\label{eqn:MansourintroTrib}
|S_n(123, 132, 3214)| = |S_n(132, 213, 1234)| = T_{n+1} \hspace{\eqnsep} (n \ge 1).
\end{equation}
Here $T_n$ is the $n$th Tribonacci number, defined by $T_0 = 0$, $T_1 = T_2 = 1$, and $T_n = T_{n-1} + T_{n-2} + T_{n-3}$ for $n \ge 3$.
In this paper we first generalize (\ref{eqn:SimionSchmidtintro}) and (\ref{eqn:MansourintroTrib}) by showing that for all $n$ and all $k \ge 2$,
\begin{equation}
\label{eqn:usintro1}
|S_n(12\ldots k, 132, 213)| = F_{k-1,n+1}
\end{equation}
and
\begin{equation}
\label{eqn:usintro2}
|S_n(123, 132, k-1\ k-2\ \ldots\ 1 k)| = F_{k-1,n+1}.
\end{equation}
Here $F_{k,n}$ is the $k$-generalized Fibonacci number defined by setting $F_{k,n} = 0$ for $n \le 0$, $F_{k,1} = 1$, and $F_{k,n} = \sum_{i=1}^k F_{k,n-i}$ for $n \ge 2$.
We give bijective proofs of (\ref{eqn:usintro1}) and (\ref{eqn:usintro2}).
We then generalize (\ref{eqn:MansourintroFib2}) by showing that for all $k \ge 4$, the set $S_n(132, k\ k-1\ \ldots 4213, 2341)$ is enumerated by a linear combination of two Fibonacci numbers plus a polynomial of degree $k - 2$.
We give explicit enumerations for $k = 4$, $k = 5$, and $k = 6$.
Here it is the $k-3$rd differences for the sequence $|S_n(132, k\ k-1\ \ldots 4213, 2341)|$ which satisfy the Fibonacci recurrence.
We go on to generalize (\ref{eqn:MansourintroFib2}) in another direction by showing that for all $k \ge 2$ and all $n \ge 1$,
$$|S_n(132, 213, 23\ldots k 1)| = \sum_{i=1}^n F_{k-2,i}.$$
In this case the first differences for the sequence $|S_n(132, 213, 23 \ldots k 1)|$ satisfy the $k-2$-generalized Fibonacci recurrence.
We conclude with some open problems and conjectures, which include conjectured generalizations of (\ref{eqn:MansourintroFib1}).

\section{Background and Notation}

In this paper, a {\it permutation} of $[n] = \{1, 2, \ldots, n\}$ is a sequence in which each element of $[n]$ appears exactly once.
We write $S_n$ to denote the set of permutations of $[n]$.
We say a permutation $\pi$ {\it avoids} a permutation $\sigma$ whenever $\pi$ contains no subsequence with all of the same pairwise comparisons as $\sigma$.
For example, the permutation 23481756 avoids 4213, but it has 4156 as a subsequence so it does not avoid 2134.
We make this idea precise in the following definition.

\begin{definition}
For any permutation $\pi \in S_n$ and any $i$ $(1 \le i \le n)$, we write $\pi(i)$ to denote the element of $\pi$ in position $i$.
We say a permutation $\pi \in S_n$ {\upshape avoids} a permutation $\sigma \in S_k$ whenever there is no sequence $1 \le i_{\sigma(1)} < i_{\sigma(2)} < \cdots < i_{\sigma(k)} \le n$ such that $\pi(i_1) < \pi(i_2) < \cdots < \pi(i_k)$.
\end{definition}

\noindent
If $\pi$ avoids $\sigma$ then $\pi$ is sometimes called a {\it restricted permutation} or a {\it pattern-avoiding permutation} and $\sigma$ is sometimes called a {\it forbidden subsequence}.
In this paper we will be interested in permutations which avoid several patterns, so for any set $R$ of permutations we write $S_n(R)$ to denote the elements of $S_n$ which avoid every element of $R$.
For any set $R$ of permutations we take $S_n(R)$ to be the empty set whenever $n < 0$ and we take $S_0(R)$ to be the set containing only the empty permutation.
When $R = \{\pi_1, \pi_2, \ldots, \pi_r\}$ we often write $S_n(R) = S_n(\pi_1, \pi_2, \dots, \pi_r)$.

For all integers $k \ge 1$, the $k$-generalized Fibonacci number $F_{k,n}$ satisfies the recurrence obtained by adding more terms to the recurrence for the Fibonacci numbers.
More specifically, we set $F_{k,n} = 0$ for all $n \le 0$ and we set $F_{k,1} = 1$.
For all $n \ge 2$ we define $F_{k,n}$ recursively by setting
\begin{equation}
\label{eqn:kgenrecurrence}
F_{k,n} = \sum_{i=1}^k F_{k,n-i} \hspace{\eqnsep} (n \ge 2).
\end{equation}
The term ``$k$-generalized Fibonacci number'' alludes to the fact that $F_{2,n} = F_n$ for all $n$.
We also observe that $F_{3,n} = T_n$, the $n$th Tribonacci number, for all $n$.
We will make use of the following combinatorial interpretation of $F_{k,n}$.

\begin{proposition}
\label{prop:kgencombinatorial}
The number of tilings of a $1 \times n$ rectangle with tiles of size $1 \times 1$, $1 \times 2$, \ldots, $1 \times k$ is the $k$-generalized Fibonacci number $F_{k,n+1}$.
\end{proposition}
\begin{proof}
The result is immediate for $n \le 1$, so it suffices to show that the number of such tilings satisfies (\ref{eqn:kgenrecurrence}).
To do this, observe there is a one-to-one correspondence between tilings of a $1 \times (n - i)$ rectangle and tilings of a $1 \times n$ rectangle in which the right-most tile has length $i$.
Therefore, if we count tilings of a $1 \times n$ rectangle according to the length of the right-most tile, we find the number of such tilings satisfies (\ref{eqn:kgenrecurrence}), as desired.
\end{proof}

For more information concerning the $k$-generalized Fibonacci numbers, see \cite{Flores}, \cite{Gabai}, \cite{Lynch}, \cite{Miles}, and \cite{Miller}.

\section{Two Families of Restricted Permutations Counted by the $k$-Generalized Fibonacci Numbers}

In \cite[Prop. 15]{SimionSchmidt} Simion and Schmidt show that for all integers $n$,
\begin{equation}
\label{eqn:SimionSchmidt}
|S_n(123,132,213)| = F_{n+1}.
\end{equation}
Generalizing this result, Mansour shows in \cite[Ex. 1, 3]{Mansour33k} that for all integers $n$,
\begin{displaymath}
|S_n(1234, 132, 213)| = |S_n(123, 132, 3214)| = T_{n+1}.
\end{displaymath}
Here $T_n$ is the Tribonacci number, defined by $T_0 = 0, T_1 = T_2 = 1$, and $T_n = T_{n-1} + T_{n-2} + T_{n-3}$ for $n \ge 3$.
In this section we generalize these results by lengthening one of the subsequences which is forbidden in (\ref{eqn:SimionSchmidt}).
There are three natural ways to do this:  replace 123 with $12\ldots k$, replace 132 with $1 k\ k-1\ \ldots 2$, or replace $213$ with $k-1\ k-2\ \ldots 1 k$.
With this in mind, we wish to enumerate $S_n(12\ldots k, 132, 213)$, $S_n(123, 1 k\ k-1\ \ldots 2, 213)$, and $S_n(123, 132, k-1\ k-2\ \ldots 1 k)$.

Our first step is to observe that enumerating $S_n(123, 1 k\ k-1\ \ldots 2, 213)$ is equivalent to enumerating $S_n(123, 132, k-1\ k-2\ \ldots 1 k)$.
To see this, suppose $\pi \in S_n$ and let $\pi^{rc}$ denote the reverse complement of $\pi$, which is the permutation of $[n]$ which satisfies
$$\pi^{rc}(i) = n + 1 - \pi(n + 1 - i) \hspace{\eqnsep} (1 \le i \le n).$$
It is routine to check that the map which takes $\pi$ to $\pi^{rc}$ restricts to a bijection between $S_n(123, 1 k\ k-1\ \ldots 2, 213)$ and $S_n(123, 132, k-1\ k-2\ \ldots 1 k)$ for all $n \ge 0$.
Therefore $|S_n(123, 1 k\ k-1\ \ldots 2, 213)| = |S_n(123, 132, k-1\ k-2\ \ldots k\ 1)|$ for all $n$.

We now turn our attention to the enumeration of $S_n(12\ldots k, 132, 213)$.
We begin by considering the form of a permutation in this set.

\begin{lemma}
\label{lem:12kform}
Fix $n \ge 1$, fix $k \ge 2$, and suppose $\pi \in S_n(12\ldots k, 132, 213)$.
Then $\pi$ has the form
\begin{equation}
\label{eqn:12kform}
n - \pi^{-1}(n) + 1,\ n - \pi^{-1}(n) + 2,\ldots, n, \sigma
\end{equation}
for some $\sigma \in S_{n-\pi^{-1}(n)}(12\ldots k, 132, 213)$.
\end{lemma}
\begin{proof}
Since $\pi$ avoids 213, the elements to the left of $n$ form an increasing subsequence.
Since $\pi$ avoids 132, the elements to the left of $n$ are all greater than every element to the right of $n$.
Combining these two observations gives the result.
\end{proof}

We now consider the permutations in $S_n(12\ldots k, 132, 213)$ according to the position of the largest element.
We observe that for $\pi \in S_n$ this position is given by $\pi^{-1}(n)$.

\begin{lemma}
\label{lem:12kpositionbound}
For all $n \ge 1$ and all $k \ge 2$ there are no permutations in $S_n(12\ldots k, 132, 213)$ which satisfy $\pi^{-1}(n) \ge k$.
\end{lemma}
\begin{proof}
If $\pi^{-1}(n) \ge k$ then by Lemma \ref{lem:12kform} the first $k$ elements of $\pi$ form a subsequence of type $12\ldots k$.
This contradicts the fact that $\pi$ avoids $12\ldots k$.
\end{proof}

\begin{lemma}
\label{lem:12kinversemap}
Fix $n \ge 1$, fix $k \ge 2$, and fix $i$, $(1 \le i \le k-1)$.
For all $\sigma \in S_{n-i}(12\ldots k, 132, \\ 213)$, the permutation
\begin{equation}
\label{eqn:12kpi}
n-i+1, n-i+2, \ldots, n, \sigma
\end{equation}
is in $S_n(12\ldots k, 132, 213)$.
\end{lemma}
\begin{proof}
Suppose $\sigma \in S_{n-i}(12\ldots k, 132, 213)$ and let $\pi$ denote the permutation in (\ref{eqn:12kpi});  we show $\pi \in S_n(12\ldots k, 132, 213)$.
Clearly $\pi \in S_n$, so it is sufficient to show $\pi$ avoids 132, 213, and $12\ldots k$.
Beginning with 132, suppose by way of contradiction that $\pi$ contains a subsequence of type 132.
Since $\sigma$ avoids 132, any 132 pattern in $\pi$ must use one of $n-i+1, \ldots, n$ as the 1.
There is no decrease in the subsequence $n-i+1, \ldots, n$, so some element of $\sigma$ must play the role of the 2.
But this contradicts the fact that every element of $\sigma$ is smaller than $n - i + 1$.
Arguing in a similar fashion, it is routine to show $\pi$ avoids 213 and $12\ldots k$.
\end{proof}

\begin{lemma}
\label{lem:12kbijection}
For all $n \ge 1$, all $k \ge 2$, and all $i$ $(1 \le i \le k-1)$,
the map from \\ $S_{n-i}(12\ldots k, 132, 213)$ to $S_n(12\ldots k, 132, 213)$ given by
\begin{equation}
\label{eqn:12kbijection}
\sigma \mapsto n-i+1,\ n-i+2,\ldots, n, \sigma
\end{equation}
is a bijection between $S_{n-i}(12\ldots k, 132, 213)$ and the elements of $S_n(12\ldots k, 132, 213)$ which satisfy $\pi^{-1}(n) = i$.
\end{lemma}
\begin{proof}
The map in (\ref{eqn:12kbijection}) is clearly injective.
Its range is contained in the correct set by Lemma \ref{lem:12kinversemap}.
It is surjective by Lemma \ref{lem:12kform}.
\end{proof}

Using Lemmas \ref{lem:12kpositionbound} and \ref{lem:12kbijection}, we obtain our enumeration of $S_n(12\ldots k, 132, 213)$.

\begin{theorem}
\label{thm:12kenumeration}
For all integers $n$ and all $k \ge 2$,
\begin{equation}
\label{eqn:12kenumeration}
|S_n(12\ldots k, 132, 213)| = F_{k-1,n+1}.
\end{equation}
Moreover, the generating function
$$f(x) = \sum_{n=0}^\infty |S_n(12\ldots k, 132, 213)| x^n$$
is given by
\begin{equation}
\label{eqn:12kgf}
f(x) = \frac{1}{1 - x - x^2 - \cdots - x^{k-1}}.
\end{equation}
\end{theorem}
\begin{proof}
Fix $k \ge 2$.
The result is immediate for $n \le 2$ so we assume $n \ge 3$.
Counting the elements of $S_n(12\ldots k, 132, 213)$ according to the value of $\pi^{-1}(n)$ and using Lemmas \ref{lem:12kpositionbound} and \ref{lem:12kbijection}, we find 
$$|S_n(12\ldots k, 132, 213)| = \sum_{i=1}^{k-1} |S_{n-i}(12\ldots k, 132, 213)|.$$
Comparing this with (\ref{eqn:kgenrecurrence}) we find that $F_{k-1,n+1}$ and $|S_n(12\ldots k, 132, 213)|$ satisfy the same recurrence.
Since they also have the same initial conditions, they are equal.
It is routine using (\ref{eqn:12kenumeration}) to obtain (\ref{eqn:12kgf}).
\end{proof}

We can also prove Theorem \ref{thm:12kenumeration} bijectively.

\begin{theorem}
\label{thm:12kbijection}
For all $n \ge 1$ and all $k \ge 2$ there exists a constructive bijection between $S_n(12\ldots k, 132, 213)$ and the set of tilings of a $1 \times n$ rectangle with tiles of size $1 \times 1$, $1 \times 2$, \ldots, $1 \times k-1$.
\end{theorem}
\begin{proof}
Suppose we are given such a tiling.
We construct its corresponding permutation as follows.
Fill the right-most tile with the numbers $1,2,\ldots$ from left to right, up to the length of the tile.
Fill the tile immediately to the left of the right-most tile from left to right, beginning with the smallest available number.
Repeat this process until every tile is filled.
It is routine using Lemma \ref{lem:12kbijection} to verify that this map is invertible, and that the permutation constructed avoids 132, 213, and $12\ldots k$, as desired.
\end{proof}

\begin{example}
Fix $k = 4$ and $n = 9$.
Under the bijection given in the proof of Theorem \ref{thm:12kbijection}, the permutation 978652341 corresponds to the tiling $\Box\ \Box\hspace{-2.5pt}\Box\ \Box\ \Box\ \Box\hspace{-2.5pt}\Box\hspace{-2.5pt}\Box\ \Box$ and the permutation 896754123 corresponds to the tiling $\Box\hspace{-2.5pt}\Box\ \Box\hspace{-2.5pt}\Box\ \Box\ \Box\ \Box\hspace{-2.5pt}\Box\hspace{-2.5pt}\Box$.
\end{example}

Having obtained our enumeration of $S_n(12\ldots k, 132, 213)$, we now turn our attention to the enumeration of $S_n(123, 132, k-1\ k-2\ \ldots 1 k)$.
We begin by considering the form of a permutation in this set.

\begin{lemma}
\label{lem:kk1form}
Fix $n \ge 1$, fix $k \ge 2$, and suppose $\pi \in S_n(123, 132, k-1\ k-2\ \ldots 1 k)$.
Then $\pi$ has the form
$$n-1, n-2, \ldots, n - \pi^{-1}(n) + 1, n, \sigma$$
for some $\sigma \in S_{n-\pi^{-1}(n)}(123, 132, k-1\ k-2\ \ldots 1 k)$.
\end{lemma}
\begin{proof}
Since $\pi$ avoids $123$, the elements to the left of $n$ form a decreasing subsequence.
Since $\pi$ avoids 132, the elements to the left of $n$ are all greater than every element to the right of $n$.
Combining these two observations gives the result.
\end{proof}

We now consider permutations in $S_n(123, 132, k-1\ k-2\ \ldots\ 1 k)$ according to the position of the largest element.
We observe that for $\pi \in S_n$ this position is given by $\pi^{-1}(n)$.

\begin{lemma}
\label{lem:kk1positionbound}
For all $n \ge 1$ and all $k \ge 2$ there are no permutations in $S_n(123, 132, k-1\ k-2\ \ldots\ 1 k)$ which satisfy $\pi^{-1}(n) \ge k$.
\end{lemma}
\begin{proof}
If $\pi^{-1}(n) \ge k$ then by Lemma \ref{lem:kk1form} the first $k-1$ elements of $\pi$ combine with $n$ to form a subsequence of type $k-1\ k-2\ \ldots 1 k$.
This contradicts the fact that $\pi$ avoids $k-1\ k-2\ \ldots 1 k$.
\end{proof}

\begin{lemma}
\label{lem:kk1inversemap}
Fix $n \ge 1$, fix $k \ge 2$, and fix $i$ $(1 \le i \le k-1)$.
For all $\sigma \in S_{n-i}(123, 132, k-1\ k-2\ \ldots\ 1 k)$, the permutation
\begin{equation}
\label{eqn:kk1inversemap}
n-1, n-2, \ldots, n-i+1, n, \sigma
\end{equation}
is in $S_n(123, 132, k-1\ k-2\ \ldots\ 1 k)$.
\end{lemma}
\begin{proof}
Suppose $\sigma \in S_{n-i}(123, 132, k-1\ k-2\ \ldots\ 1 k)$ and let $\pi$ denote the permutation in (\ref{eqn:kk1inversemap});  we show $\pi \in S_n(123, 132, k-1\ k-2\ \ldots\ 1 k)$.
Clearly $\pi \in S_n$, so it is sufficient to show $\pi$ avoids 123, 132, and $k-1\ k-2\ \ldots\ 1 k$.
Beginning with $k-1\ k-2\ \ldots\ 1 k$, suppose by way of contradiction that $\pi$ contains a subsequence of type $k-1\ k-2\ \ldots\ 1 k$.
Since $\sigma$ avoids $k-1\ k-2\ \ldots\ 1 k$, any pattern of this type in $\pi$ must use one of $n-i+1, \ldots, n$ as the $k-1$.
But this contradicts the fact that $i \le k-1$.
Arguing in a similar fashion, it is routine to show $\pi$ avoids $123$ and $132$.
\end{proof}

\begin{lemma}
\label{lem:kk1bijection}
For all $n \ge 1$, all $k \ge 2$, and all $i$ $(1 \le i \le k-1)$, the map from $S_{n-i}(123, 132, k-1\ k-2\ \ldots\ 1 k)$ to $S_n(123, 132, k-1\ k-2\ \ldots\ 1 k)$ given by
\begin{equation}
\label{eqn:kk1bijection}
\sigma \mapsto n-1, n-2, \ldots, n-i+1, n, \sigma
\end{equation}
is a bijection between $S_{n-i}(123, 132, k-1\ k-2\ \ldots\ 1 k)$ and the elements of $S_n(123, 132, k-1\ k-2\ \ldots\ 1 k)$ which satisfy $\pi^{-1}(n) = i$.
\end{lemma}
\begin{proof}
The map in (\ref{eqn:kk1bijection}) is clearly injective.
Its range is contained in the correct set by Lemma \ref{lem:kk1inversemap}.
It is surjective by Lemma \ref{lem:kk1form}.
\end{proof}

Using Lemmas \ref{lem:kk1positionbound} and \ref{lem:kk1bijection}, we immediately obtain our enumeration of $S_n(123, 132, k-1\ k-2\ \ldots\ 1 k)$.

\begin{theorem}
\label{thm:kk1enumeration}
For all integers $n$ and all $k \ge 2$,
\begin{equation}
\label{eqn:kk1enumeration}
|S_n(123, 132, k-1\ k-2\ \ldots\ 1 k)| = F_{k-1,n+1}.
\end{equation}
Moreover, the generating function
$$f(x) = \sum_{n=0}^\infty |S_n(123, 132, k-1\ k-2\ \ldots 1 k)| x^n$$
is given by
\begin{equation}
\label{eqn:kk1gf}
f(x) = \frac{1}{1 - x - \cdots - x^{k-1}}.
\end{equation}
\end{theorem}
\begin{proof}
Fix $k \ge 2$.
The result is immediate for $n \le 2$ so we assume $n \ge 3$.
Counting the elements of $S_n(123, 132, k-1\ k-2\ \ldots\ 1 k)$ according to the value of $\pi^{-1}(n)$ and using Lemmas \ref{lem:kk1positionbound} and \ref{lem:kk1bijection}, we find
$$|S_n(123, 132, k-1\ k-2\ \ldots\ 1 k)| = \sum_{i=1}^{k-1} |S_{n-i}(123, 132, k-1\ k-2\ \ldots\ 1 k)|.$$
Comparing this with (\ref{eqn:kgenrecurrence}) we find that $F_{k-1,n+1}$ and $|S_n(123, 132, k-1\ k-2\ \ldots\ 1 k)|$ satisfy the same recurrence.
Since they also have the same initial conditions, they are equal.
It is routine using (\ref{eqn:kk1enumeration}) to obtain (\ref{eqn:kk1gf}).
\end{proof}

We can also prove Theorem \ref{thm:kk1enumeration} bijectively.

\begin{theorem}
\label{thm:kk1bijection}
For all $n \ge 1$ and all $k \ge 2$ there exists a constructive bijection between $S_n(123, 132, k-1\ k-2\ \ldots 1 k)$ and the set of tilings of a $1 \times n$ rectangle with tiles of size $1 \times 1$, $1 \times 2$, \ldots, $1 \times k-1$.
\end{theorem}
\begin{proof}
Suppose we are given such a tiling.
We construct its corresponding permutation as follows.
Let $m$ denote the length of the right-most tile.
Fill this tile with the numbers $1, 2, \ldots, m$ from left to right in the order $m-1, m-2, \ldots, 1, m$.
Fill the tile immediately to the left of the right-most tile in the same way, using the smallest available numbers.
Repeat this process until every tile is filled.
It is routine using Lemma \ref{lem:kk1bijection} to verify that this map is invertible, and that the permutation constructed avoids 123, 132, and $k-1\ k-2\ \ldots 1 k$.
\end{proof}

\begin{example}
Fix $k = 5$ and $n = 9$.
Under the bijection given in the proof of Theorem \ref{thm:kk1bijection}, the permutation 879643251 corresponds to the tiling $\Box\hspace{-2.5pt}\Box\hspace{-2.5pt}\Box\ \Box\ \Box\hspace{-2.5pt}\Box\hspace{-2.5pt}\Box\hspace{-2.5pt}\Box\ \Box$ and the permutation 986754213 corresponds to the tiling $\Box\ \Box\ \Box\hspace{-2.5pt}\Box\ \Box\ \Box\ \Box\hspace{-2.5pt}\Box\hspace{-2.5pt}\Box$.
\end{example}

\section{A Family of Restricted Permutations Counted by a Fibonacci Number Plus a Polynomial}
\label{sec:fibpolyfamily}

In \cite[Ex. 3]{Mansour33k} Mansour shows that
$$|S_n(132, 213, 2341)| = F_{n+2} - 1 \hspace{\eqnsep} (n \ge 1).$$
In this section we generalize this result by extending the permutation 213 to permutations of the form $k\ k-1\ \ldots 4213$.
We show that for any $k \ge 4$, the set $S_n(132, k\ k-1\ \ldots 4213, 2341)$ is enumerated by a linear combination of two Fibonacci numbers plus a polynomial of degree $k-2$.
We give explicit enumerations for $k = 4$, $k = 5$, and $k = 6$.
We begin by considering the form of a permutation in $S_n(132, k\ k-1\ \ldots 4213, 2341)$.

\begin{lemma}
\label{lem:kk-1form}
Fix $n \ge 1$, fix $k \ge 4$, and suppose $\pi \in S_n(132, k\ k-1\ \ldots 4213, 2341)$.
Then the following hold.
\renewcommand\labelenumi{{\upshape (\roman{enumi}) }}
\begin{enumerate}
\item If $\pi^{-1}(n) = 1$ and $n > 1$ then $\pi$ has the form
\begin{equation}
\label{eqn:kk-1form1}
n, \sigma
\end{equation}
for some $\sigma \in S_{n-1}(132, k-1\ k-2\ \ldots 4213, 2341)$.
\item If $\pi^{-1}(n) = 2$ and $n > 2$ then $\pi$ has the form
\begin{equation}
\label{eqn:kk-1form2}
n-1,\ n, \sigma
\end{equation}
for some $\sigma \in S_{n-2}(132, k-1\ k-2\ \ldots 4213, 2341)$.
\item If $3 \le \pi^{-1}(n) \le k-3$ and $\pi^{-1}(n) < n$ then $\pi$ has the form
\begin{equation}
\label{eqn:kk-1form3k-3}
n-1,\ n-2, \ldots,n-\pi^{-1}(n)+1,\ n, \sigma
\end{equation}
for some $\sigma \in S_{n-\pi^{-1}(n)}(132, k-\pi^{-1}(n)+1\ \ldots 4213, 2341)$.
\item
If $k-2 \le \pi^{-1}(n) \le n-1$ then $\pi$ has the form
\begin{equation}
\label{eqn:kk-1formk-2n-1}
n-1,\ n-2,\ldots,n-\pi^{-1}(n)+1,\ n,\sigma
\end{equation}
for some $\sigma \in S_{n-\pi^{-1}(n)}(132, 213, 2341)$.
\item
If $\pi^{-1}(n) = n$ then $\pi$ has the form
\begin{equation}
\label{eqn:kk-1formn}
\sigma, n
\end{equation}
for some $\sigma \in S_{n-1}(132, k\ k-1\ \ldots 4213, 2341)$.
\end{enumerate}
\end{lemma}
\begin{proof}
(i)
Observe it is sufficient to show $\sigma$ avoids $k-1\ k-2\ \ldots 4213$.
Suppose by way of contradiction that $\sigma$ contains a pattern of this type.
Combining this pattern with $n$ gives a pattern of type $k\ k-1\ \ldots 4213$ in $\pi$, contradicting the fact that $\pi$ avoids $k\ k-1\ \ldots 4213$.
Therefore $\sigma \in S_{n-1}(132, k-1\ k-2\ \ldots 4213, 2341)$, as desired.

(ii)
Since $\pi$ avoids 132, the element to the left of $n$ is greater than every element to the right of $n$, so it must be $n-1$.
Therefore $\pi$ has the form in (\ref{eqn:kk-1form2}) for some $\sigma \in S_{n-2}$.
Arguing as in the proof of (i), we now find $\sigma \in S_{n-2}(132, k-1\ k-2\ \ldots 4213, 2341)$.

(iii)
Since $\pi$ avoids 132, the elements to the left of $n$ are all greater than every element to the right of $n$.
Since $\pi$ avoids 2341 and $\pi^{-1}(n) < n$, the elements to the left of $n$ form a decreasing subsequence.
Therefore $\pi$ has the form in (\ref{eqn:kk-1form3k-3}) for some $\sigma \in S_{n - \pi^{-1}(n)}$.
To show $\sigma \in S_{n-\pi^{-1}(n)}(132, k-\pi^{-1}(n) + 1\ \ldots 4213, 2341)$ it is sufficient to show $\pi$ avoids $k - \pi^{-1}(n) + 1\ \ldots 4213$.
Suppose by way of contradiction that $\sigma$ contains a pattern of this type.
Combining this pattern with $n-1,\ n-2, \ldots, n-\pi^{-1}(n)+1$ gives a pattern of type $k\ k-1\ \ldots 4213$ in $\pi$, contradicting the fact that $\pi$ avoids $k\ k-1\ \ldots 4213$.
Therefore $\sigma \in S_{n-\pi^{-1}(n)}(132, k-\pi^{-1}(n) + 1\ \ldots 4213, 2341)$.

(iv)
This is similar to the proof of (iii).

(v)
This is immediate.
\end{proof}

Digressing for a moment, we now use Lemma \ref{lem:kk-1form} to give a bijective proof of Mansour's result.

\begin{proposition}
\label{prop:Mansourfn+2-1}
For all $n \ge 1$, there exists a constructive bijection between $S_n(132, 213, \\ 2341)$ and the set of tilings of a $1 \times (n+1)$ rectangle with tiles of size $1 \times 1$ and $1 \times 2$ using at least one $1 \times 2$ tile.
In particular, for all $n \ge 1$,
$$|S_n(132, 213, 2341)| = F_{n+2} - 1.$$
\end{proposition}
\begin{proof}
Suppose we are given such a tiling.
We construct its corresponding permutation as follows.
Replace the right-most $1 \times 2$ tile with a 1, and fill the (necessarily $1 \times 1$) tiles to the right of the 1 with $2, 3, \ldots$ from left to right.
Now fill the right-most empty tile with the smallest numbers available, placing them in the tile from left to right in increasing order.
Repeat this process until every tile is filled.
It is routine using Lemma \ref{lem:kk-1form} to verify that this map is invertible, and that the permutation constructed avoids 132, 213, and 2341, as desired.
\end{proof}

\begin{example}
Under the bijection given in the proof of Proposition \ref{prop:Mansourfn+2-1}, the permutation 87564123 corresponds to the tiling $\Box\ \Box\ \Box\hspace{-2.5pt}\Box\ \Box\ \Box\hspace{-2.5pt}\Box\ \Box\ \Box$ and the permutation 86745321 corresponds to the tiling $\Box\ \Box\hspace{-2.5pt}\Box\ \Box\hspace{-2.5pt}\Box\ \Box\ \Box\ \Box\hspace{-2.5pt}\Box$.
\end{example}

We now return to our study of the permutations in $S_n(132, k\ k-1\ \ldots 4213, 2341)$.

\begin{lemma}
\label{lem:kk-1inverse}
The following hold for all $n \ge 1$ and all $k \ge 4$.
\renewcommand\labelenumi{{\upshape (\roman{enumi}) }}
\begin{enumerate}
\item
For all $\sigma \in S_{n-1}(132, k-1\ k-2\ \ldots 4213, 2341)$, the permutation
\begin{equation}
\label{eqn:kk-1inverse1}
n, \sigma
\end{equation}
is in $S_n(132, k\ k-1\ \ldots 4213, 2341)$.
\item
For all $\sigma \in S_{n-2}(132, k-1\ k-2\ \ldots 4213, 2341)$, the permutation
\begin{equation}
\label{eqn:kk-1inverse2}
n-1,\ n, \sigma
\end{equation}
is in $S_n(132, k\ k-1\ \ldots 4213, 2341)$.
\item
Suppose $3 \le i \le k-3$ and $i < n$.
Then for all $\sigma \in S_{n-i}(132, k-i+1\ \ldots 4213, 2341)$, the permutation
\begin{equation}
\label{eqn:kk-1inverse3k-3}
n-1,\ n-2,\ \ldots,n-i+1,\ n, \sigma
\end{equation}
is in $S_n(132, k\ k-1\ \ldots 4213, 2341)$.
\item
Suppose $k-2 \le i \le n - 1$.
Then for all $\sigma \in S_{n-i}(132, 213, 2341)$, the permutation
\begin{equation}
\label{eqn:kk-1inversek-2n-1}
n-1,\ n-2,\ldots,n-i+1,\ n,\sigma
\end{equation}
is in $S_n(132, k\ k-1\ \ldots 4213, 2341)$.
\item
For all $\sigma \in S_{n-1}(132, k\ k-1\ \ldots 4213, 2341)$, the permutation
\begin{equation}
\label{eqn:kk-1inversen}
\sigma, n
\end{equation}
is in $S_n(132, k\ k-1\ \ldots 4213, 2341)$.
\end{enumerate}
\end{lemma}
\begin{proof}
This is similar to the proof of Lemma \ref{lem:kk1inversemap}.
\end{proof}

\begin{lemma}
\label{lem:kk-1bijection}
The following hold for all $n \ge 1$ and all $k \ge 4$.
\renewcommand\labelenumi{{\upshape (\roman{enumi}) }}
\begin{enumerate}
\item
The map from $S_{n-1}(132, k-1\ k-2\ \ldots 4213, 2341)$ to $S_n(132, k\ k-1\ \ldots 4213, 2341)$ given by
\begin{equation}
\label{eqn:kk-1bijection1}
\sigma \mapsto n, \sigma
\end{equation}
is a bijection between $S_{n-1}(132, k-1\ k-2\ \ldots 4213, 2341)$ and the elements of $S_n(132, k\ k-1\ \ldots 4213, 2341)$ which satisfy $\pi^{-1}(n) = 1$.
\item
The map from $S_{n-2}(132, k-1\ k-2\ \ldots 4213, 2341)$ to $S_n(132, k\ k-1\ \ldots 4213, 2341)$ given by
\begin{equation}
\label{eqn:kk-1bijection2}
\sigma \mapsto n-1,\ n, \sigma
\end{equation}
is a bijection between $S_{n-2}(132, k-1\ k-2\ \ldots 4213, 2341)$ and the elements of $S_n(132, k\ k-1\ \ldots 4213, 2341)$ which satisfy $\pi^{-1}(n) = 2$.
\item
Suppose $3 \le i \le k-3$ and $i \le n - 1$.
Then the map from $S_{n-i}(132, k-i+1\ \ldots 4213, 2341)$ to $S_n(132, k\ k-1\ \ldots 4213, 2341)$ given by
\begin{equation}
\label{eqn:kk-1bijection3k-3}
\sigma \mapsto n-1,\ n-2,\ \ldots,n-i+1,\ n, \sigma
\end{equation}
is a bijection between $S_{n-i}(132, k-i+1\ \ldots 4213, 2341)$ and the elements of $S_n(132, k\ k-1\ \ldots 4213, 2341)$ which satisfy $\pi^{-1}(n) = i$.
\item
Suppose $k-2 \le i < n$.
Then the map from $S_{n-i}(132, 213, 2341)$ to $S_n(132, k\ k-1\ \ldots 4213, 2341)$ given by
\begin{equation}
\label{eqn:kk-1bijectionk-2n-1}
\sigma \mapsto n-1,\ n-2,\ \ldots,n-i+1,\ n, \sigma
\end{equation}
is a bijection between $S_{n-i}(132, 213, 2341)$ and the elements of $S_n(132, k\ k-1\ \ldots 4213, \\ 2341)$ which satisfy $\pi^{-1}(n) = i$.
\item
The map from $S_{n-1}(132, k\ k-1\ \ldots 4213, 2341)$ to $S_n(132, k\ k-1\ \ldots 4213, 2341)$ given by
\begin{equation}
\label{eqn:kk-1bijectionn}
\sigma \mapsto \sigma, n
\end{equation}
is a bijection between $S_{n-1}(132, k\ k-1\ \ldots 4213, 2341)$ and the elements of $S_n(132, k\ k-1\ \ldots 4213, 2341)$ which satisfy $\pi^{-1}(n) = n$.
\end{enumerate}
\end{lemma}
\begin{proof}
These maps are clearly injective.
Their ranges are contained in the correct sets by Lemma \ref{lem:kk-1inverse}.
They are surjective by Lemma \ref{lem:kk-1form}.
\end{proof}

Using the last lemma, we obtain a recurrence for $|S_n(132, k\ k-1\ \ldots 4213, 2341)|$.

\begin{proposition}
\label{prop:kk-1recurrence}
For all $n \ge 1$ and all $k \ge 4$,
\begin{eqnarray*}
\lefteqn{|S_n(132, k\ k-1\ \ldots 4213, 2341)| = }\\[2ex]
& & |S_{n-1}(132, k-1\ k-2\ \ldots 4213, 2341)| + |S_{n-2}(132, k-1\ k-2\ \ldots 4213, 2341)| \\[2ex] 
& & + \sum_{i=3}^{min(k-3, n-1)} |S_{n-i}(132, k-i+1\ \ldots 4213, 2341)| + \sum_{i=k-2}^{n-1} |S_{n-i}(132, 213, 2341)| \\[2ex]
& & + |S_{n-1}(132, k\ k-1\ \ldots 4213, 2341)|. \\
\end{eqnarray*}
\end{proposition}
\begin{proof}
Count the elements of $S_n(132, k\ k-1\ \ldots 4213, 2341)$ according to the value of $\pi^{-1}(n)$ and use Lemma \ref{lem:kk-1bijection}.
\end{proof}

Using Proposition \ref{prop:kk-1recurrence} we can now enumerate $S_n(132, 4213, 2341)$.
To do so, we will need the following results concerning various sums of Fibonacci numbers, generating functions for certain subsequences of the Fibonacci numbers, and generating functions for certain sequences of binomial coefficients.
We omit the proofs, which are routine.

\begin{proposition}
\label{prop:FibSums}
For all $k \ge 1$ and all $n \ge k$,
\begin{equation}
\label{eqn:FibSums}
\sum_{i=k}^n F_i = F_{n+2} - F_{k+1}.
\end{equation}
\end{proposition}

\begin{proposition}
\label{prop:Fibgens}
For all $k \ge 2$, the generating function 
$$f(x) = \sum_{n=0}^\infty F_{n+k} x^n$$
is given by
$$f(x) = \frac{F_{k-1} x + F_k}{1 - x - x^2}.$$
\end{proposition}

\begin{proposition}
\label{prop:binomialgens}
For all $k \ge 0$, the generating function
$$f(x) = \sum_{n=0}^\infty {{n+1} \choose {k}} x^n$$
is given by
$$f(x) = \frac{x^{max(k-1, 0)}}{(1-x)^{k+1}}$$.
\end{proposition}

We now enumerate $S_n(132, 4213, 2341)$.

\begin{theorem}
\label{thm:Mansourk=4}
For all $n \ge 1$,
$$|S_n(132, 4213, 2341)| = F_{n+5} - {{n+1} \choose {2}} - 2 {{n+1} \choose {1}} - 2.$$
Moreover, the generating function
$$f(x) = \sum_{n=0}^\infty |S_n(132, 4213, 2341)| x^n$$
is given by
\begin{equation}
\label{eqn:Mansourk=4gf}
f(x) = \frac{x^4 - x^3 - 3x^2 + 3x - 1}{(x^2 + x - 1) (1 - x)^3}.
\end{equation}
\end{theorem}
\begin{proof}
The result is immediate for $n = 1,2,$ and $3$ so we assume $n \ge 4$.
Setting $k = 4$ in Proposition \ref{prop:kk-1recurrence} and simplifying the result, we find
$$|S_n(132, 4213, 2341)| = |S_{n-1}(132, 4213, 2341)| + \sum_{i=1}^{n-1} |S_i(132, 213, 2341)|.$$
Using Proposition \ref{prop:Mansourfn+2-1} and (\ref{eqn:FibSums}) to simplify the right side we find
\begin{eqnarray*}
|S_n(132, 4213, 2341)| &=& |S_{n-1}(132, 4213, 2341)| + \sum_{i=1}^{n-1} (F_{i+2} - 1) \\
&=& |S_{n-1}(132, 4213, 2341)| + F_{n+3} - n - 2.
\end{eqnarray*}
Applying the last line to itself repeatedly and using (\ref{eqn:FibSums}) to simplify the result, we find
\begin{eqnarray*}
|S_n(132, 4213, 2341)| &=& |S_3(132, 4213, 2341)| + \sum_{i=4}^n (F_{i+3} - i - 2) \\
&=& F_{n+5} - {{n+1} \choose {2}} - 2 {{n+1} \choose {1}} - 2,
\end{eqnarray*}
as desired.
It is routine using this and Propositions \ref{prop:Fibgens} and \ref{prop:binomialgens} to obtain (\ref{eqn:Mansourk=4gf}).
\end{proof}

Using the same approach, we now enumerate $S_n(132, 54213, 2341)$.

\begin{theorem}
\label{thm:Mansourk=5}
For all $n \ge 2$,
$$|S_n(132, 54213, 2341)| = 3 F_{n+5} - 2 {{n+1} \choose {3}} - 4 {{n+1} \choose {2}} - 3 {{n+1} \choose {1}} - 14.$$
Moreover, the generating function
$$f(x) = \sum_{n=0}^\infty |S_n(132, 54213, 2341)| x^n$$
is given by
\begin{equation}
\label{eqn:Mansourk=5gf}
f(x) = \frac{x^7 - 3x^5 + x^4 + 2x^3 - 6x^2 + 4x-1}{(x^2 + x - 1)(1-x)^4}.
\end{equation}
\end{theorem}
\begin{proof}
The result is immediate for $n = 2,3$, and $4$ so we assume $n \ge 5$.
Setting $k = 5$ in Proposition \ref{prop:kk-1recurrence} and simplifying the result, we find
\begin{eqnarray*}
\lefteqn{|S_n(132, 54213, 2341)| = } \\[2ex]
& & |S_{n-1}(132, 54213, 2341)| + |S_{n-1}(132, 4213, 2341)| \\[2ex]
& & + |S_{n-2}(132, 4213, 2341)| + \sum_{i=3}^{n-1} |S_{n-i}(132, 213, 2341)|.\\
\end{eqnarray*}
Using Proposition \ref{prop:Mansourfn+2-1}, Theorem \ref{thm:Mansourk=4}, and (\ref{eqn:FibSums}) to simplify the terms on the right side we find
$$|S_n(132, 54213, 2341)| = |S_{n-1}(132, 54213, 2341)| + F_{n+5} + F_{n+1} - n^2 - 3n - 3.$$
Applying this to itself repeatedly and using (\ref{eqn:FibSums}) again, we find
\begin{eqnarray*}
|S_n(132, 54213, 2341)| &=& |S_4(132, 54213, 2341)| + \sum_{i=5}^n \left( F_{i+5} + F_{i+1} - i^2 - 3i - 3\right) \\
&=& 3 F_{n+5} - 2 {{n+1} \choose {3}} - 4 {{n+1} \choose {2}} - 3 {{n+1} \choose {1}} - 14, \\
\end{eqnarray*}
as desired.
It is routine using this and Propositions \ref{prop:Fibgens} and \ref{prop:binomialgens} to obtain (\ref{eqn:Mansourk=5gf}).
\end{proof}

Finally, we enumerate $S_n(132, 654213, 2341)$.

\begin{theorem}
\label{thm:Mansourk=6}
For all $n \ge 2$,
$$|S_n(132, 654213, 2341)| = 5 F_{n+6} + F_{n+4} - 4 {{n+1} \choose {4}} - 7 {{n+1} \choose {3}} - 5 {{n+1} \choose {2}} - 27 {{n+1} \choose {1}} - 8.$$
Moreover, the generating function
$$f(x) = \sum_{n=0}^\infty |S_n(132, 654213, 2341)| x^n$$
is given by
\begin{equation}
\label{eqn:Mansourk=6gf}
f(x) = \frac{2x^8 - x^7 - x^6- 5 x^5-x^4+8x^3 - 10x^2 + 5x - 1}{(x^2 + x - 1)(1 - x)^5}.
\end{equation}
\end{theorem}
\begin{proof}
The result is immediate for $n = 2,3$, and $4$ so we assume $n \ge 5$.
Setting $k = 6$ in Proposition \ref{prop:kk-1recurrence} and simplifying the result, we find
\begin{eqnarray*}
\lefteqn{|S_n(132, 654213, 2341)| = } \\[2ex]
& & |S_{n-1}(132, 654213, 2341)| + |S_{n-1}(132, 54213, 2341)| + |S_{n-2}(132, 54213, 2341)| \\[2ex]
& & + |S_{n-3}(132, 4213, 2341)| + \sum_{i=4}^{n-1} |S_{n-i}(132, 213, 2341)| \\
\end{eqnarray*}
Using Proposition \ref{prop:Mansourfn+2-1}, Theorems \ref{thm:Mansourk=4} and \ref{thm:Mansourk=5}, and (\ref{eqn:FibSums}) to simplify the terms on the right side we find
$$|S_n(132, 654213, 2341)| = |S_{n-1}(132, 654213, 2341)| + 3 F_{n+5} + F_{n+2} + F_n - \frac{2}{3} n^3 - \frac{3}{2} n^2 - \frac{17}{6} n - 27.$$
Applying this to itself repeatedly and using (\ref{eqn:FibSums}) again, we find
\begin{eqnarray*}
\lefteqn{|S_n(132, 654213, 2341)|} \\
&=& |S_4(132, 654213, 2341)| + \sum_{i=5}^n \left( 3 F_{i+5} + F_{i+2} + F_i - \frac{2}{3} i^3 - \frac{3}{2} i^2 - \frac{17}{6} i - 27 \right) \\
&=& 5 F_{n+6} + F_{n+4} - 4 {{n+1} \choose {4}} - 7 {{n+1} \choose {3}} - 5 {{n+1} \choose {2}} - 27 {{n+1} \choose {1}} - 8, \\
\end{eqnarray*}
as desired.
It is routine using this and Propositions \ref{prop:Fibgens} and \ref{prop:binomialgens} to obtain (\ref{eqn:Mansourk=6gf}).
\end{proof}

Using the same approach, one can (with enough patience and computational power) enumerate $S_n(132, k\ k-1\ \ldots\ 4213, 2341)$ for any $k \ge 4$.
Rather than give explicit enumerations for larger values of $k$, we content ourselves with the following result concerning the form these enumerations will take.

\begin{theorem}
\label{thm:FibPolyFamily}
For all $n$ and $k$ such that $n \ge k \ge 4$, the set $S_n(132, k\ k-1\ \ldots 4213, 2341)$ is enumerated by a linear combination of two Fibonacci numbers plus a polynomial of degree $k - 2$.
\end{theorem}
\begin{proof}
Assume $n \ge k$; we argue by induction on $k$.
The cases $k=4$, $k=5$, and $k=6$ are done in Theorems \ref{thm:Mansourk=4}, \ref{thm:Mansourk=5}, and \ref{thm:Mansourk=6}, respectively.
Therefore we assume $k \ge 7$ and that the result holds for all $k' < k$.
Inspecting the recurrence in Proposition \ref{prop:kk-1recurrence}, we find by induction that with the exception of $|S_{n-1}(132, k\ k-1\ \ldots 4213, 2341)|$, each term on the right side may be written as a linear combination of two Fibonacci numbers plus a polynomial of degree at most $k - 3$.
In particular, we find
\begin{eqnarray*}
\lefteqn{|S_n(132, k\ k-1\ \ldots 4213, 2341)| =} \\[2ex]
& & |S_{n-1}(132, k\ k-1\ \ldots 4213, 2341)| + p_{k-3}(n) + a F_{n+r} + b F_{n+s}
\end{eqnarray*}
for some polynomial $p_{k-3}(n)$ of degree $k - 3$ and constants $a$, $b$, $r$, and $s$.
Applying this to itself repeatedly we find
\begin{eqnarray*}
\lefteqn{|S_n(132, k\ k-1\ \ldots 4213, 2341)| =} \\[2ex]
& & |S_4(132, k\ k-1\ \ldots 4213, 2341)| + \sum_{i=5}^n \left( p_{k-3}(i) + a F_{i+r} + b F_{i+s} \right). \\
\end{eqnarray*}
Using (\ref{eqn:FibSums}) to simplify the right side of this equation, we find $|S_n(132, k\ k-1\ \ldots 4213, 2341)|$ may be written as a linear combination of two Fibonacci numbers plus a polynomial of degree $k - 2$, as desired.
\end{proof}

We can use Theorem \ref{thm:FibPolyFamily} to give a simple description of the connection between the sequence $|S_n(132, k\ k-1\ \ldots 4213, 2341)|$ and the Fibonacci numbers.
For a given sequence $a_0, a_1, a_2, \ldots$ we refer to the sequence $a_1 - a_0, a_2 - a_1, a_3 - a_2, \ldots$ as the {\it first difference sequence} for $a_0, a_1, \ldots$.
Similarly, we refer to the sequence obtained from $a_0, a_1, \ldots$ by taking differences $i$ times as the {\it $i$th difference sequence} for $a_0, a_1, \ldots$.
Although the exact enumeration of $S_n(132, k\ k-1\ 4213, 2341)$ seems to be quite complicated, we have the following result concerning the $k-3$rd differences for the sequence $|S_n(132, k\ k-1\ 4213, 2341)|$.

\begin{proposition}
\label{prop:k-3differences}
For all $k \ge 4$, the $k-3$rd differences for $|S_n(132, k\ k-1\ \ldots 4213, 2341)|$ satisfy the Fibonacci recurrence.
\end{proposition}
\begin{proof}
Observe that if $a_0, a_1, \ldots$ satisfies the Fibonacci recurrence then so do all of its difference sequences.
Also observe that if $a_i$ is a polynomial of degree $n$ in $i$ then the $n+1$st difference sequence for $a_0, a_1, \ldots$ consists entirely of zeroes.
The result follows by combining these two observations with Theorem \ref{thm:FibPolyFamily}.
\end{proof}

\section{A Family of Restricted Permutations Counted by Sums of $k$-Generalized Fibonacci Numbers}
\label{sec:kgendiff}

In \cite[Ex. 3]{Mansour33k} Mansour shows that
$$|S_n(132, 213, 2341)| = F_{n+2} - 1 \hspace{\eqnsep} (n \ge 1).$$
In this section we generalize this result by lengthening the permutation 2341 to permutations of the form $23\ldots k 1$.
We show that for all $k \ge 4$,
$$|S_n(132, 213, 23\ldots k 1)| = \sum_{i=1}^n F_{k-2,i} \hspace{\eqnsep} (n \ge 1).$$
This complements Proposition \ref{prop:k-3differences} in that it implies that the first differences of the sequence $|S_n(132, 213, 23\ldots k 1)|$ are the $k-2$-generalized Fibonacci numbers.
We begin by considering the form of a permutation in $S_n(132, 213, 23\ldots k1)$.

\begin{lemma}
\label{lem:23k1form}
Fix $n \ge 1$, fix $k \ge 2$, and suppose $\pi \in S_n(132, 213, 23\ldots k 1)$.
Then $\pi$ has the form
\begin{equation}
\label{eqn:23k1form}
n - \pi^{-1}(n) + 1, n - \pi^{-1}(n) + 2, \ldots, n, \sigma
\end{equation}
for some $\sigma \in S_{n - \pi^{-1}(n)}(132, 213, 23\ldots k 1)$.
\end{lemma}
\begin{proof}
This is similar to the proof of Lemma \ref{lem:12kform}.
\end{proof}

We now consider the permutations in $S_n(132, 213, 23\ldots k1)$ according to the position of the largest element.
We observe that for $\pi \in S_n$ this position is given by $\pi^{-1}(n)$.

\begin{lemma}
\label{lem:23k1positionbound}
Fix $n \ge 1$, fix $k \ge 2$, and suppose $\pi \in S_n(132, 213, 23\ldots k 1)$.
If $\pi^{-1}(n) \ge k-1$ then $\pi = 12\ldots n$.
\end{lemma}
\begin{proof}
Suppose $\pi \in S_n(132, 213, 23\ldots k 1)$ and $\pi^{-1}(n) \ge k-1$.
If $\pi^{-1}(n) < n$ then by Lemma \ref{lem:23k1form} any $k-2$ entries of $\pi$ to the left of $n$, combined with $n$ and any entry to the right of $n$, form a pattern of type $23\ldots k 1$.
Therefore $\pi^{-1}(n) = n$ and $\pi = 12\ldots n$ by the form of (\ref{eqn:23k1form}).
\end{proof}

\begin{lemma}
\label{lem:23k1inversemap}
Fix $n \ge 1$, fix $k \ge 2$, and fix $i$ $(1 \le i \le min(k-2, n-1))$.
For all $\sigma \in S_{n-i}(132, 213, 23\ldots k 1)$, the permutation
\begin{equation}
\label{eqn:23k1inversemap}
n-i+1, n-i+2, \ldots, n, \sigma
\end{equation}
is in $S_n(132, 213, 23\ldots k 1)$.
\end{lemma}
\begin{proof}
Suppose $\sigma \in S_{n-i}(132, 213, 23\ldots k 1)$ and let $\pi$ denote the permutation in (\ref{eqn:23k1inversemap});  we show $\pi \in S_n(132, 213, 23\ldots k 1)$.
Clearly $\pi \in S_n$, so it is sufficient to show $\pi$ avoids 132, 213, and $23\ldots k 1$.
Beginning with $23\ldots k 1$, suppose by way of contradiction that $\pi$ contains a subsequence of type $23\ldots k 1$.
Since $\sigma$ avoids $23\ldots k 1$, any pattern of this type in $\pi$ must use one of $n-i+1, \ldots, n$ as the 2.
Since $i \le k-1$, this pattern must also use an element of $\sigma$ as the $k$, contradicting the form of (\ref{eqn:23k1form}).
Arguing in a similar fashion, it is routine to show $\pi$ avoids 132 and 213.
\end{proof}

\begin{lemma}
\label{lem:23k1bijection}
For all $n \ge 1$, all $k \ge 2$, and all $i$ $(1 \le i \le min(k-2, n-1))$, the map from $S_{n-i}(132, 213, 23\ldots k 1)$ to $S_n(132, 213, 23\ldots k 1)$ given by
\begin{equation}
\label{eqn:23k1bijection}
\sigma \mapsto n-i+1, n-i+2, \ldots, n, \sigma
\end{equation}
is a bijection between $S_{n-i}(132, 213, 23\ldots k 1)$ and the elements of $S_n(132, 213, 23\ldots k 1)$ which satisfy $\pi^{-1}(n) = i$.
\end{lemma}
\begin{proof}
The map in (\ref{eqn:23k1bijection}) is clearly injective.
Its range is contained in the correct set by Lemma \ref{lem:23k1inversemap}.
It is surjective by Lemma \ref{lem:23k1form}.
\end{proof}

Using the Lemmas \ref{lem:23k1positionbound} and \ref{lem:23k1bijection} we obtain a recurrence for $|S_n(132, 213, 23\ldots k1)|$.

\begin{proposition}
\label{prop:fn2-1}
Fix $k \ge 3$.
For all $n \ge k-1$ we have
\begin{equation}
\label{eqn:fibrecurrenceplusone}
|S_n(132, 213, 23\ldots k 1)| = 1 + \sum_{i=1}^{k-2} |S_{n-i}(132, 213, 23\ldots k 1)|.
\end{equation}
\end{proposition}
\begin{proof}
Count the elements of $S_n(132, 213, 23\ldots k 1)$ according to the value of $\pi^{-1}(n)$ and use Lemmas \ref{lem:23k1positionbound} and \ref{lem:23k1bijection}.
\end{proof}

Using (\ref{eqn:fibrecurrenceplusone}) we now show that the first differences for the sequence $|S_n(132, 213, 23\ldots k1)|$ are the $k-2$-generalized Fibonacci numbers.
To do this we will need the following result of Simion and Schmidt.

\begin{proposition}
\label{prop:SS2n-1}
\cite[Prop. 8]{SimionSchmidt}
For all $n \ge 1$,
\begin{equation}
\label{eqn:SS2n-1}
|S_n(132, 213)| = 2^{n-1}.
\end{equation}
\end{proposition}

\begin{proposition}
For all $k \ge 3$ and all $n \ge 1$,
\begin{equation}
\label{eqn:kgendifference}
|S_n(132, 213, 23\ldots k1)| - |S_{n-1}(132, 213, 23\ldots k1)| = F_{k-2,n}.
\end{equation}
\end{proposition}
\begin{proof}
Suppose $1 < n < k$.
In this case $|S_n(132, 213, 23\ldots k1)| = |S_n(132, 213)|$, so by (\ref{eqn:SS2n-1}) we have $|S_n(132, 213, 23\ldots k1)| - |S_{n-1}(132, 213, 23\ldots k1)| = 2^{n-2} = F_{k-2,n}$.
Now suppose $n \ge k$.
In this case by (\ref{eqn:fibrecurrenceplusone}) we have
$$1 = |S_n(132, 213, 23\ldots k 1)| - \sum_{i=1}^{k-2} |S_{n-i}(132, 213, 23\ldots k 1)|$$
and 
$$1 = |S_{n-1}(132, 213, 23\ldots k 1)| - \sum_{i=1}^{k-2} |S_{n-i-1}(132, 213, 23\ldots k 1)|.$$
Subtract the second line from the first to find that
\begin{eqnarray*}
\lefteqn{|S_n(132, 213, 23\ldots k 1)| - |S_{n-1}(132, 213, 23\ldots k 1)| =} \\
& & \sum_{i=1}^{k-2} \left( |S_{n-i}(132, 213, 23\ldots k 1)| - |S_{n-i-1}(132, 213, 23\ldots k 1)| \right).
\end{eqnarray*}
Comparing this with (\ref{eqn:kgenrecurrence}) we find that $|S_n(132, 213, 23\ldots k 1)| - |S_{n-1}(132, 213, 23\ldots k 1)|$ and $F_{k-2,n}$ satisfy the same recurrence.
Since they also have the same initial conditions, they are equal.
\end{proof}

\begin{proposition}
For all $k \ge 2$ and all $n \ge 1$,
\begin{equation}
\label{eqn:secondtolast}
|S_n(132, 213, 23\ldots k1)| = \sum_{i=1}^n F_{k-2,i}.
\end{equation}
Moreover, the generating function
$$f(x) = \sum_{n=0}^\infty |S_n(132, 213, 23\ldots k1)| x^n$$
is given by
\begin{equation}
\label{eqn:last}
f(x) = \frac{x}{(1-x)(1 - x - x^2 - \ldots - x^{k-2})} + 1.
\end{equation}
\end{proposition}
\begin{proof}
By (\ref{eqn:kgendifference}) we have
$$|S_n(132, 213, 23\ldots k1)| = |S_{n-1}(132, 213, 23\ldots k1)| + F_{k-2,n}.$$
Applying this to itself repeatedly gives (\ref{eqn:secondtolast}).
It is now routine using (\ref{eqn:fibrecurrenceplusone}) and Proposition \ref{prop:SS2n-1} to obtain (\ref{eqn:last}).
\end{proof}

\section{Open Problems and Conjectures}

We conclude with some open problems and conjectures suggested by the results in this paper.

\begin{enumerate}
\item In \cite[Ex. 1]{Mansour33k} Mansour shows that for all $n \ge 1$,
$$|S_n(123, 132, 4213)| = F_{n+2} - 1.$$
Based on numerical evidence, we conjecture that the following also hold.
\begin{equation}
\label{eqn:op1}
|S_n(123, 2431, 4213)| = F_{n+3} + 2n - 9 \hspace{\eqnsep} (n \ge 3)
\end{equation}
\begin{equation}
\label{eqn:op2}
|S_n(123, 35421, 4213)| = F_{n+3} + F_{n+1} + 2 {{n+1} \choose {2}} - n - 23 \hspace{\eqnsep} (n \ge 5)
\end{equation}
One ought to be able to prove these conjectures using generating trees.
(See \cite{WestGenTrees} for more information on generating trees.)

\item
Generalize (\ref{eqn:op1}) and (\ref{eqn:op2}) to an infinite family $R_3, R_4, \ldots$ of sets of forbidden subsequences such that $|S_n(R_i)|$ is given by a linear combination of two Fibonacci numbers plus a polynomial for all $i \ge 3$.

\item
In section \ref{sec:fibpolyfamily} we show that the $k-3$rd difference sequence for the sequence $|S_n(132, k\ k-1\ \ldots 4213, 2341)|$ satisfies the Fibonacci recurrence.
In section \ref{sec:kgendiff} we show that the first difference sequence for $|S_n(132, 213, 23\ldots k1)|$ satisfies the $k-2$-generalized Fibonacci recurrence.
For all $i \ge 2$ and all $k \ge 3$, find a set $R_{i,k}$ such that the $i$th difference sequence for $|S_n(R_{i,k})|$ satisfies the $k$-generalized Fibonacci recurrence. 

\item
Give bijective proofs of Theorems \ref{thm:Mansourk=4}, \ref{thm:Mansourk=5}, and \ref{thm:Mansourk=6}.

\item
Give an explicit enumeration of $S_n(132, k\ k-1\ \ldots 4213, 2341)$ for all $k \ge 4$.
\end{enumerate}

\bigskip

\begin{Large}
\noindent
{\bf Acknowledgement}
\end{Large}

\bigskip
The author would like to thank Darla Kremer for many helpful comments, suggestions, and conversations.

\bibliographystyle{jalg}
\bibliography{references}

\end{document}